\documentclass[preprint,12pt]{elsarticle}
\usepackage{amsthm}
\newtheorem{definition}{\textbf{Definition}}

\newtheorem*{meta*}{Meta-Conjecture}
\newtheorem{theorem}{\textbf{Theorem}}
\newtheorem*{question*}{\textbf{Question}}
\newtheorem{corollary}{\textbf{Corollary}}
\newtheorem{remark}{\textbf{Remark}}

\newtheorem{lemma}{\textbf{Lemma}}

\usepackage{mathrsfs}

\def\a {\alpha}
\def\o {\omega}

\def\N {\mathbb{N}}
\def\Z {\mathbb{Z}}
\def\Q {\mathbb{Q}}
\def\L {\mathbb{L}}
\def\R {\mathbb{R}}
\def\U {\mathbb{U}_m}

\def\QQ {\overline{\Q}}
\def\Ls {\mathbb{L}_{\scriptsize{strong}}}

\def\cL {\mathcal{L}}
\def\e {\epsilon}

\usepackage{amssymb}
\usepackage{amsmath}
\usepackage{ragged2e}
\usepackage[symbol]{footmisc}

\journal{...}

\begin{document}

\begin{frontmatter}

\title{On the Arithmetic Behavior of Liouville Numbers under Rational Maps}

\author[ap]{Ana Paula Chaves}{\ead{apchaves@ufg.br}}
\author[dm]{Diego Marques\corref{CA}}{\ead{diego@mat.unb.br}}
\author[pt]{Pavel Trojovsk\'y}{\ead{pavel.trojovsky@uhk.cz}}
\address[ap]{Instituto de Matem\'atica e Estat\'istica, Universidade Federal de Goi\'as, Brazil}
\address[dm]{Departamento de Matem\'atica, Universidade de Bras\'ilia, Brazil}
\address[pt]{Faculty of Science, University of Hradec Kr\'alov\'e, Czech Republic}

\cortext[CA]{Corresponding author}

\begin{abstract}
In 1972, Alnia\c{c}ik proved that every strong Liouville number is mapped into the set of $U_m$-numbers, for any non-constant rational function with coefficients belonging to an $m$-degree number field. In this paper, we generalize this result by providing a larger class of Liouville numbers (which, in particular, contains the strong Liouville numbers) with this same property (this set is sharp is a certain sense).

\end{abstract}


\end{frontmatter}

\section{Introduction}

The beginning of the transcendental number theory happened in 1844, when Liouville \cite{liouville} showed that algebraic numbers are not ``well-approximated"\ by rationals. More precisely, if $\alpha$ is an $n$-degree real algebraic number (with $n>1$), then there exists a positive constant $C$, depending only on $\a$, such that $|\alpha -p/q| > Cq^{-n}$, for all rational number $p/q$. By using this result, Liouville was able to explicit, for the first time, examples of transcendental numbers (the now called {\it Liouville numbers}). 

A real number $\xi$ is called a Liouville number if there exist infinitely many rational numbers $(p_n/q_n)_n$, with $q_n> 1$, such that
\begin{center}
$\displaystyle 0 < \left| \xi - \frac{p_n}{q_n} \right| < \frac{1}{q_n^{\omega_n}}$,
\end{center}
for some sequence of real numbers $(\o_n)_n$ which tends to $\infty$ as $n\to \infty$. The set of the Liouville numbers is denoted by $\mathbb{L}$.
 
The first example of a Liouville number (and consequently, of a transcendental number) is the so-called {\it Liouville constant} defined by the convergent series $\ell = \sum_{n\geq 1}10^{-n!}$ (i.e., the decimal with $1$'s in each factorial position and $0$'s otherwise). In 1962, Erd\"os \cite{erdos} proved that every real number can be written as the sum of two Liouville numbers (indeed, besides being a completely topological property of $G_{\delta}$-dense sets, Erd\"{o}s was able to provide an explicit proof based only on the definition of Liouville numbers).

Let $\xi$ be a real irrational number and $(p_k/q_k)_k$ be the sequence of the convergents of its continued fraction. It is well-known that the definition of Liouville numbers is equivalent to: for every positive integer $n$, there exist infinitely many positive integers $k$ such that $q_{k+1}  > q_k^n$. Thus, we have the following subclass of Liouville numbers.

\begin{definition}
Let $\xi$ be a real irrational number and $(p_k/q_k)_k$ be the sequence of the convergents of its continued fraction. The number $\xi$ is said to be a {\emph strong Liouville number} if, for every $n$, there exists $N$ (depending on $n$) such that $q_{k+1}>q_k^n$ for all $k > N$. For convenience, we shall denote by $\Ls$ the set of all strong Liouville numbers. 
\end{definition}

Observe that the previous definition is equivalent to say that $\xi\in \Ls$ if for the sequence $(\omega_k)_k$ defined by
\[
\left| \xi - \frac{p_k}{q_k} \right|=q_k^{-\o_k},
\]
it holds that $\lim\inf_{k\to \infty}\o_k=\infty$.

Erd\"os also raised the problem whether every real number can be written as the sum of two strong Liouville numbers. This question was answered by Petruska \cite{petruska}, who proved that the sum of two strong Liouville numbers is either a rational or a Liouville number. In the same paper, he also showed that for any integer $a\geq 2$ and any sequence of positive integers $(v_n)_n$ such that $v_{n+1} \geq 2v_n+2$, the number 
\begin{equation}\label{petruskanumber}
\displaystyle\sum_{n \geq 1}a^{-v_n} 
\end{equation}
is not a strong Liouville number. In particular, the Liouville constant is not a strong Liouville number (take $a=10$ and $v_n=n!$). 

A related open problem is to describe the additive and multiplicative group generated by the strong Liouville numbers. For instance, $\xi^2$ is never a strong Liouville number if $\xi\in \Ls$.

In 1932, Mahler \cite{mahler} split the set of transcendental numbers into three disjoint sets: $S$-, $T$- and $U$-numbers (according to their approximation by  algebraic numbers). In particular, $U$-numbers generalize the concept of Liouville numbers. 

Let $\omega^*_n(\xi)$ be the supremum of the positive real numbers $\omega^*$ for which there exist infinitely many real $n$-degree algebraic numbers $\alpha$ satisfying
\[
0 < |\xi - \alpha| < H(\alpha)^{-\omega^*-1},
\]
where $H(\alpha)$ (so-called the {\it height} of $\alpha$) is the maximum of the absolute values of the coefficients of the minimal polynomial of $\a$ (over $\Z$). If $\o^*_m(\xi) = \infty$ and $\o^*_n(\xi) < \infty$, for all $1 \leq n <m$, the number $\xi$ is said to be a $U^*_m$-number. Actually, this is the definition of Koksma's $U^*_m$-numbers (see \cite{koksma}). However, it is well known that the sets of $U_m$- and $U_m^*$-numbers are the same \cite{bugeaud} (let us denote the set of $U_m$-numbers by $\U$). We still remark that the set of Liouville numbers is precisely the set of $U_1$-numbers.

The existence of $U_m$-numbers, for all $m\geq 1$, was proved by LeVeque \cite{leveque}. In fact, he proved that $\sqrt[m]{(3+\ell)/4}$ is a $U_m$-number, for all $m\geq 1$. In 2014, Chaves and Marques \cite{chavesmarques} explicit $U_m$-numbers in a more natural way: as the product of certain $m$-degree algebraic numbers by the Liouville constant $\ell$ (for example, $\sqrt[m]{2/3}\cdot \ell$ is a $U_m$-number, for all $m\geq 1$). In 1972, Alnia\c{c}ik \citep{alniacik}  proved that a much more general fact holds for strong Liouville numbers. Indeed, if $\xi$ is a strong Liouville number, then $F(\xi)$ is a $U_m$-number, for all rational function $F(x)$ with coefficients in an $m$-degree number field.

In this paper, we generalize this Alnia\c{c}ik's result to a larger class of Liouville numbers which, in particular, contains all strong Liouville numbers. More precisely

\begin{theorem} \label{main}
Let $(\o_n)_n$ be a sequence of positive real numbers which tends to $\infty$ as $n\to \infty$. Let $\xi$ be a real number such that there exists an infinite sequence of rational numbers $(p_n/q_n)_n$, satisfying
\begin{equation} \label{ineq1}
\left| \xi - \frac{p_n}{q_n} \right| < H\left(\frac{p_n}{q_n}\right)^{-\omega_n},
\end{equation}
where $H(p_{n+1}/q_{n+1}) \leq H(p_{n}/q_{n})^{O(\omega_n)}$, for all $n\geq 1$. Then, for any irreducible rational function $F(x)\in K(x)$, with $[K:\Q]=m$, the number $F(\xi)$ is a $U_m$-number.
\end{theorem}

From now on, the set of the Liouville numbers satisfying the conditions of the theorem will be denoted by $\mathcal{L}$. Note also that the previous theorem implies that the set $\mathcal{L}$ is not $G_{\delta}$-dense (in $\R$) and so it is a small set in the Hausdorff sense (in fact, on the contrary, i.e., if $\cL$ is a $G_{\delta}$-dense set, then every real number can be written as sum of two elements of $\cL$. However, for instance, if $\sqrt{2}=\ell_1+\ell_2$, for some $\ell_1, \ell_2\in \cL$, we would obtain an absurdity as $\mathbb{U}_2\ni \sqrt{2}-\ell_1=\ell_2\in \mathbb{L}$).

From the definition, we can deduce that $\Ls\subseteq \cL$ and that $\cL\backslash \Ls$ is an uncountable set. Indeed, the second assertion follows from the fact that all numbers in the form (\ref{petruskanumber}) with $v_{n+1}/v_n\to \infty$ as $n\to \infty$ belongs to $\cL\backslash \Ls$ (the details of this fact will be provided in the proof of the first corollary). For the first assertion, if $\xi$ is a strong Liouville number, and $(\o_k)_k$ is defined by $|\xi-p_k/q_k|=q_k^{-\o_k}$ (where $(p_k/q_k)_k$ are the convergents of the continued fraction of $\xi$), then
\[
\dfrac{1}{q_k^{2\o_k}}<q_k^{-\o_k}=\left| \xi - \frac{p_k}{q_k} \right|<\min\left\{\dfrac{1}{q_kq_{k+1}},\dfrac{1}{q_k^{\o_k/2}}\right\}.
\]
Thus $|\xi-p_k/q_k|<q_k^{-\o_k/2}$ and $q_{k+1}<q_k^{O(\o_k/2)}$ which implies that $\xi\in \cL$ as desired. 

\begin{remark}
We observe that Theorem \ref{main} is sharp is a certain sense (see the last section). In fact, the conclusion of the theorem is not valid if we instead consider Liouville numbers as in that statement, but with the weaker condition
\[
H(p_{n+1}/q_{n+1}) \leq H(p_{n}/q_{n})^{O((\omega_n)^{1+\epsilon})}
\]
for any previously fixed $\epsilon>0$.
\end{remark}

Now, we shall list some immediate consequences of Theorem \ref{main}.

\begin{corollary}\label{c1}
If $\xi$ is a Liouville number of the form \eqref{petruskanumber}, where $v_{n+1}/v_n$ tends to infinity as $n\to \infty$ (in particular, $\xi\in \mathcal{L}$), then  $F(\xi)$ is a $U_m$-number, for any non-constant rational function $F(x)$ with coefficients in an $m$-degree number field. 
\end{corollary}

The next result is related to a discussion of the second author with Bugeaud in a personal communication.

\begin{corollary}\label{c2}
We have that
\begin{itemize}
\item[(i)] There exist infinitely many $U_m$-numbers $\xi$ such that for any positive divisor $d$ of $m$, there is a positive integer $k$ such that $\xi^k$  is a $U_d$-number.
\item[(ii)] There exist infinitely many $U_m$-numbers $\xi$ such that $\xi^k$ is a $U_m$-number, for all $k\geq 1$.
\end{itemize}
\end{corollary}

\begin{corollary}\label{c3}
For any integer $m\geq 1$, let $\psi: \U\times \QQ(x)\to \N$ be the function defined of the following form: For $\xi\in \U$ and $F(x)\in \QQ(x)$, the number $\psi(\xi,F)$ is defined as the type of the $U$-number $F(\xi)$. Then $\psi^{-1}(s)$ is an infinite set, for any positive integer $s$ (in particular, $\psi$ is a surjective function). Indeed, for any $s\geq 1$, the series
\[
\displaystyle\sum_{\xi\in \pi_1(\psi^{-1}(s))}\dfrac{1}{\xi}
\]
diverges. Here, as usual, $\pi_1(x,y)=x$ is a projection function.
\end{corollary}

Let us describe in a few words the main idea of the proof of Theorem \ref{main}. First, we shall use some results on heights to show the well-approximation of $F(\xi)$ by $m$-degree algebraic numbers (which are, in fact, the image of the rational approximants of $\xi$ by $F$). This implies that $F(\xi)$ is a $U$-number of type at most $m$. To show that its type is exactly $m$, we shall use some asymptotic estimates together with a ``Liouville-type"\ inequality due to Bombieri (related to the distance between distinct algebraic numbers) to deduce that $|F(\xi)-\gamma|>H(\gamma)^{\nu}$ for all (but finitely many) $n$-degree algebraic numbers (for any $1\leq n<m$). In conclusion, $F(\xi)$ must be a $U_m$-number.

\section{Auxiliary results}\label{sec2}

In this section, we shall present some facts which will be essential ingredients in the proof of our theorem. The first one is due to \' I\c{c}en \cite{íçen} which provides an inequality relating the height of algebraic numbers.

\begin{lemma}\label{l1}
Let $\alpha_1,\ldots,\alpha_k$ be algebraic numbers belonging to a $g$-degree number field $K$, and let $F(y,x_1,\ldots,x_k)$ be a polynomial with integer coefficients and with degree $d>1$ in the variable $y$. If $\eta$ is an algebraic number such that $F(\eta,\alpha_1,\ldots,\alpha_k)=0$, then the degree of $\eta$ is at most $dg$, and 
\begin{equation} \label{içendes}
H(\eta) \leq 3^{2dg + (l_1+\cdots + l_k)g} \cdot H^g \cdot H(\a_1)^{l_1g} \cdots H(\a_k)^{l_kg}, 
\end{equation}
where $H$ is the maximum of the absolute values of the coefficients of $F$ (the height of $F$), $l_i$ is the degree of $F$ in the variable $x_i$, for $i=1,\ldots, k$.
\end{lemma}

The next lemma is a {\it Liouville-type} inequality due to Bombieri \cite{bombieri}.

\begin{lemma}\label{l2}
Let $\a_1, \a_2$ be distinct algebraic numbers of degrees $n_1$ and $n_2$, respectively. Then
\begin{equation}\label{algeb}
|\a_1-\a_2| > (4n_1n_2)^{-3n_1n_2}H(\a_1)^{-n_2}H(\a_2)^{-n_1}.
\end{equation}
\end{lemma}

Also, we shall use the relations below (which can be proved by using some inequalities in \cite[Section 3.4]{waldschmidt}. See also \cite{guguActa}). 

\begin{lemma}\label{l3}
For $\alpha_1,\ldots, \alpha_n$ belonging to some $d$-degree number field, the following inequalities hold
\begin{enumerate}
\item[{(a)}] $H(\alpha_1+\cdots+\alpha_n) \leq \left[ 2n(d+1)^{\frac{n}{2}} \right]^d\cdot H(\alpha_1)^d \cdots H(\alpha_n)^d$;
\item[(b)] $H(\alpha_1\cdots\alpha_n) \leq \left[ 2(d+1)^{\frac{n}{2}} \right]^d\cdot H(\alpha_1)^d \cdots H(\alpha_n)^d$.
\end{enumerate} 
\end{lemma}

The next lemma is a result due to Alnia\c{c}ik \cite{alniacik} and will be useful to precise the exact degree of the algebraic approximants.

\begin{lemma}\label{l4}
Let $F(x)$ be an irreducible rational function with coefficients in an $m$-degree number field $K$. Then $F(\a)$ is a primitive element of the field extension $K/\Q$, for all rational number $\a$, with finitely many exceptions. In particular, $F(\a)$ is an $m$-degree algebraic number, for all sufficiently large rational number $\a$.
\end{lemma}

Now, we are ready to deal with the proof of our result. 

\section{The proofs}

\subsection{The proof of Theorem \ref{main}}

In what follows, the implied constants in $\ll$ and $\gg$ depend only on $m, k, l$ and the coefficients of the rational function $F(x)$. In particular, we shall write $k\gg 1$ in order to avoid the unnecessary repetition of ``for all sufficiently large $k$". Also, $c_1, c_2,\ldots$ will denote real constants $>1$.

We start by taking
\begin{center}
 $\a_k:=\dfrac{p_k}{q_k}$, $F(z) := \dfrac{P(z)}{Q(z)}$ and $\gamma_k:= F(\alpha_k)$.  
\end{center}
Let $C>1$ be a real number such that  $H(\alpha_{k+1}) \leq H(\alpha_{k})^{C\omega_k}$ for $k\gg 1$. Now, we can apply the Mean Value Theorem for $F(x)$ and the points $\xi$ and $\a_k$ (when both numbers belong to some very small interval which does not contain any pole or zero of $F(x)$ and $F'(x)$) to obtain
\begin{equation} \label{mvt}
|F(\xi) - F(\a_k)| \ll \left| \xi - \frac{p_k}{q_k} \right|.
\end{equation}
Now, we combine \eqref{ineq1} and \eqref{mvt} to get
 \begin{equation} \label{desig1}
 |F(\xi) - \gamma_k| \ll H(\alpha_k)^{-\omega_k}.
 \end{equation}
On the other hand, we have that (by Lemma \ref{l3})
\[
H(\gamma_k) = H(P(\a_k) Q(\a_k)^{-1}) \ll H(P(\a_k))^m H(Q(\a_k))^m.
\]
Again from Lemma \ref{l3}, we can rewrite the last inequality as
\begin{equation} \label{hgamma}
H(\gamma_k) \ll H(\a_k)^{2m^2},
\end{equation}
and so,
\begin{equation} \label{desig2}
H(\a_k)^{-\omega_k} < c_1^{\omega_k} H(\gamma_k)^{-\frac{\omega_k}{2m^2}}.
\end{equation}
From \eqref{desig1} and \eqref{desig2}, we have
\[
 |F(\xi) - \gamma_k| \ll \frac{c_1^{\omega_k}}{H(\gamma_k)^{\frac{\omega_k}{4m^2}}} \cdot H(\gamma_k)^{-\frac{\omega_k}{4m^2}}.
\]
Due to Lemma \ref{l4}, we have that $F(\a_k)$ is an $m$-degree algebraic number for $k\gg 1$. By using that $\lim \sup_{k\to \infty} H(\gamma_k)=\infty$ (since we have infinitely many $\a_k$'s), we deduce that
\begin{equation}
|F(\xi) - \gamma_k| \ll H(\gamma_k)^{-\frac{\omega_k}{4m^2}},
\end{equation}
for $k\gg 1$. This yields $\omega^*_m(F(\xi)) = \infty$ and so $F(\xi)$ is a $U$-number of type $\leq m$.

In order to prove that the type of $F(\xi)$ (as a $U$-number) is exactly $m$, we shall show that $F(\xi)$ is not well-approximated by $n$-degree algebraic numbers, for $n<m$. In fact, let $\gamma$ be an algebraic number of degree $n$ ($<m$). Then, from Lemma \ref{l2}, we have
\begin{equation} \label{alg2}
|\gamma - \gamma_k| \gg H(\gamma)^{-m} H(\gamma_k)^{-n}.
\end{equation}
Also, for some $k \gg 1$, it holds that
\begin{equation} \label{desig3}
H(\gamma_k) \leq H(\gamma)^{16Cm^6} < H(\gamma_{k+1}).
\end{equation} 
Now, by combining the inequality in \eqref{hgamma} with the fact that $H(\alpha_{k+1}) \leq H(\alpha_{k})^{C\omega_k}$, we get 
\begin{equation} \label{desig4}
H(\gamma_{k+1}) \ll H(\alpha_{k})^{2Cm^2\omega_k}.
\end{equation}

Let us rewrite $\gamma_k = P(\a_k)/Q(\a_k)=\sum_{j=0}^la_j\a_k^j/\sum_{j=0}^rb_j\a_k^j$ as
\begin{equation}\label{poli}
\gamma_kb_0+\gamma_kb_1\a_k + \cdots +\gamma_kb_r\a_k^r - a_0 -a_1\a_k - \cdots - a_l\a_k^l = 0.
\end{equation}
So, $G(\a_k,\gamma_k,b_0,\ldots,b_r,a_0,\ldots,a_l)=0$, where  $G(y,x_1,\ldots , x_{r+l+3})$ is defined by
\[
G(y,x_1,\ldots , x_{r+l+3}) = \displaystyle\sum_{j=1}^{r+1}x_1x_{j+1}y^{j-1} - \displaystyle\sum_{j=1}^{l+1}x_{r+j+2}y^{j-1}.
\]
Then, by \eqref{poli} and Lemma \ref{l1}, we arrive at
\[
H(\a_k) \ll H(\gamma_k)^m,
\]
and so
\[
H(\a_k)^{2Cm^2\omega_k} \leq c_2^{\omega_k} \cdot H(\gamma_k)^{2Cm^3\omega_k}.
\]
Thus, by combining the last inequality together with \eqref{desig3}, we obtain
\begin{equation} \label{desigprinc}
H(\gamma_k) \leq H(\gamma)^{16Cm^6} \ll c_2^{\omega_k} \cdot H(\gamma_k)^{2Cm^3\omega_k}.
\end{equation}
By a straightforward manipulation, we get
\[
H(\gamma)^{-m} \gg c_3^{-\omega_k}\cdot H(\gamma_k)^{-\frac{\omega_k}{8m^2}},
\]  
and this allows us to rewrite \eqref{alg2} as
\begin{equation} \label{desig5}
|\gamma - \gamma_k| >c_4 \cdot c_3^{-\omega_k} H(\gamma_k)^{-\frac{\omega_k}{8m^2}-n}.
\end{equation}
By taking $H(\gamma) \gg 1$ (and so $H(\a_k)\gg 1$) such that
\begin{equation} \label{desig6}
c_4 \cdot c_3^{-\omega_k} > 2H(\gamma_k)^{-\frac{\omega_k}{8m^2}+n},
\end{equation}
then \eqref{desig5} and \eqref{desig6} give
\[
|\gamma - \gamma_k|  > 2H(\gamma_k)^{-\frac{\omega_k}{4m^2}} > 2 |F(\xi) - \gamma_k|.
\]
On the other hand, we have
\begin{eqnarray*}
|F(\xi) - \gamma| & \geq &  |\gamma - \gamma_k| - |F(\xi) - \gamma_k| \\
& > & \frac{1}{2} |\gamma - \gamma_k| > \frac{c_3}{2} H(\gamma)^{-m} \end{eqnarray*}
and therefore
\[
|F(\xi) - \gamma| \gg  H(\gamma)^{-m-16Cnm^6}.
\]
In conclusion, $\omega_n^*(\xi) < m+16Cnm^6<\infty$, for all $n<m$, which implies that $F(\xi)$ is a $U_m$-number. This completes the proof.\qed

\subsection{The proof of the corollaries}

\subsubsection{Proof of the Corollary \ref{c1}} In order to prove this result, it suffices to see that if $\xi$ is a Liouville number of the form \eqref{petruskanumber}, where $v_{k+1}/v_k$ tends to infinity as $k\to \infty$, then for $q_k=a^{v_k}$, we have
\[
\left|\xi-\dfrac{p_k}{q_k}\right|\ll a^{-v_{k+1}}=q_k^{-v_{k+1}/v_k}.
\]
Also, $q_{k+1}=a^{v_{k+1}}=q_k^{v_{k+1}/v_k}$ and so the result follows from Theorem \ref{main}.\qed

\subsubsection{Proof of the Corollary \ref{c2}} To prove this corollary, we first observe that if $\ell$ is the Liouville constant, then $\ell^k$ satisfies the hypotheses of Theorem \ref{main}, for all $k\geq 1$. 
\begin{itemize}
    \item[(i).] Then, by our main result, we have that $\sqrt[m]{p}\cdot \ell$ is a $U_m$-number, for every prime number $p$. So, if $d\mid m$ ($d>0$), then $(\sqrt[m]{p}\cdot \ell)^{m/d}=\sqrt[d]{p}\cdot \ell^{m/d}$ is a $U_d$-number (again by Theorem \ref{main}).
    \item[(ii).] For a prime number $p$, take the $U_m$-number $\xi:=(1+\sqrt[m]{p})\ell$, then $\xi^k$ is still a $U_m$-number, for all $k\geq 1$, by the fact that $(1+\sqrt[m]{p})^k$ is always an $m$-degree algebraic number, together with our main theorem. \qed
\end{itemize}

\subsubsection{Proof of the Corollary \ref{c3}} First, note that the function $\psi$ is well defined. In fact, since any $F(x)\in \QQ(X)$ can be considered with coefficients belonging to an $r$-degree number field $K$. Then, if $\xi$ is a $U_m$-number with $m$-degree algebraic approximants $\a_k$, then, by the same arguments as in the proof of  Theorem \ref{main}, we conclude that almost all $F(\a_k)$ (unless $\a_k$ is a pole of $F$) will be approximants of $F(\xi)$. Since $F(\a_k)$ is an algebraic number with degree at most $rm$, then infinitely many of them will have the same degree, say $d$ ($\leq rm$). Thus, $F(\xi)$ is a $U$-number of type at most $d$ and so $\psi$ is well defined. Now, for a positive integer $s$, take $\xi=\sqrt[m]{2}\cdot \ell$ and $F(x)=\sqrt[s]{2}\cdot x^m$. Clearly, by our main theorem, we have that $\xi\in \U$ and $F(\xi)\in \mathbb{U}_s$. Moreover, for any prime number $p$, the pair $(\sqrt[m]{p}\cdot \ell,\sqrt[s]{2}\cdot x^m)\in \psi^{-1}(s)$ and so
\[
\displaystyle\sum_{\xi\in \pi_1(\psi^{-1}(s))}\dfrac{1}{\xi}\geq \displaystyle\sum_{p\ \mbox{{\scriptsize prime}}}\dfrac{1}{\sqrt[m]{p}\cdot \ell}\geq \dfrac{1}{\ell}\displaystyle\sum_{p\ \mbox{{\scriptsize prime}}}\dfrac{1}{p}=\infty.
\]
\qed

\section{Final remarks}

For $m>1$, the set of $U_m$-numbers is not a $G_{\delta}$-dense subset of $\R$, and thus it is not possible to use Erd\"{o}s' non-constructive proof to ensure the decomposition of real numbers as sum of two $U_m$-numbers (see \cite{ASTU} for results on the Borel hierarchy of the $S$, $T$ and $U$-numbers). In 1990, Alnia\c cik \cite{al} proved that every real number (except possibly Liouville numbers) can be written as sum of two $U_2$-numbers. His approach depends on the theory of continued fractions. Shortly thereafter, Pollington \cite{pol} adapted a method used by Schmidt (to ensure the existence of $T$-numbers, see \cite{schi}) and showed that for any natural number $n$, every real number may be decomposed as the sum of two $U_n$-numbers. Pollington's proof is not too easy to follow (since it depends on an inductive construction of intervals and also on some probability facts). Now, we shall provide another (possible) recipe to prove this. In fact, the Erd\"{o}s' constructive proof can be rephrased as $\R = \cL_1+\cL_1$, where $\cL_{\e}$ is the set of all real numbers $\xi$ such that there exists a sequence $(\o_n)_n$ of positive numbers (tending to $\infty$ as $n\to \infty$) and an infinite sequence of rational numbers $(p_n/q_n)_n$, satisfying
\begin{equation} 
\left| \xi - \frac{p_n}{q_n} \right| < H\left(\frac{p_n}{q_n}\right)^{-\omega_n},
\end{equation}
where $H(p_{n+1}/q_{n+1}) \leq H(p_{n}/q_{n})^{O((\omega_n)^{1+\e})}$, for all $n\geq 1$. Note that $\cL=\cL_0$.

So, in order to give an ``easier"\ (more direct) proof that $\R=\U+\U$, it would be enough to prove the following ``conjecture":
\begin{meta*}
For any $m>1$, there exists an $m$-degree algebraic number $\a$ such that $\a\cdot \cL_1\subseteq \U$.
\end{meta*}
In fact, $\R=\U+\U$ would follow directly from the fact that for any real number $x$ and for such $\a$ (as in the ``conjecture"), we have that $x/\a=\ell_1+\ell_2$, where $\ell_i\in \cL_1$ $(i=1,2)$. Thus, $x=\a\ell_1+\a\ell_2\in \U+\U$ as desired. 

Unfortunately, this ``meta-conjecture"\ is false. In fact, Erd\"{o}s' proof yields that any, nonzero, $m$-degree algebraic number ($m>1$) $\a$ can be written in the form $\ell_1\cdot \ell_2$, where $\ell_i\in \cL_1$ $(i=1,2)$. Since $\cL_1$ is closed by inversion ($x\mapsto 1/x$), we have that $\a\ell_3$ does not belong to $\U$ (since $\a\ell_3=\ell_2\in \cL_1$, for $\ell_3=1/\ell_2$).

Thus, Theorem \ref{main} is sharp in the sense in which is not possible to weaken the condition $q_{k+1}\leq q_k^{O(\o_k)}$ in order to obtain the same conclusion of that theorem. This is an immediate consequence of the previous discussion together with the fact that $\R=\cL_{\e}+\cL_{\e}$, for all $\e>0$ (it is straightforward to mimic Erd\"{o}s' construction proof to deduce this).

All these suggest the following question:
\begin{question*}
Under what conditions on the convergents (approximants) of a Liouville number, it would be possible to construct a non-$G_{\delta}$-dense set $\textsc{L}$ such that $\cL\subseteq \textsc{L}\subseteq \L$ and $F(\textsc{L})\subseteq \U$, for all non-constant rational function $F(x)$ with coefficients in a $m$-degree number field?
\end{question*}

\section*{Acknowledgement}

A part of the preparation of this paper was made during a visit of A.P.C and D.M to IMPA during its Summer Program 2019. They thank this institution for its support, hospitality and excellent working conditions. A.P.C was supported in part by CNPq Universal 01/2016 - 427722/2016-0 grant. D.M is supported by a Productivity and Research scholarship from CNPq-Brazil. P.T has been  supported by Specific Research Project of Faculty of Science, University of Hradec Kr\'alov\'e, No. 2101, 2018.

\end{document}